\documentclass[12pt]{article}
\usepackage{amsfonts}
\usepackage{a4wide}
\usepackage{amsmath,amscd,amsthm,a4,amssymb}

\begin{document}

\begin{center}
{\large \textbf{On the Diliberto-Straus algorithm for the uniform
approximation by a sum of two algebras}}

\

\textbf{Aida Kh. Asgarova}$^{1}$\textbf{\ and Vugar E. Ismailov}\footnote{%
Corresponding author}$^{2}$ \vspace{1mm}

$^{1,2}${Institute of Mathematics and Mechanics, National Academy of
Sciences of Azerbaijan, Az-1141, Baku, Azerbaijan} \vspace{1mm}

e-mail: $^{1}$aidaasgarova@gmail.com, $^{2}${vugaris@mail.ru}\vspace{5mm}
\end{center}

\textbf{Abstract.} In 1951, Diliberto and Straus \cite{4} proposed a
levelling algorithm for the uniform approximation of a bivariate function,
defined on a rectangle with sides parallel to the coordinate axes, by sums
of univariate functions. In the current paper, we consider the problem of
approximation of a continuous function defined on a compact Hausdorff space
by a sum of two closed algebras containing constants. Under reasonable
assumptions, we show the convergence of the Diliberto-Straus algorithm. For
the approximation by sums of univariate functions, it follows that
Diliberto-Straus's original result holds for a large class of compact convex
sets.

\bigskip



\bigskip

\begin{center}
{\large \textbf{1. Introduction}}
\end{center}

Let $E$ be a Banach space and $U$ and $V$ be closed subspaces thereof. In
addition, let $A$ and $B$ be proximity maps (best approximation operators)
acting from $E$ onto $U$ and $V$, respectively. We are interested in
algorithmic methods for computing the distance to a given element $z\in E$
from $U+V$. Historically, there is a procedure called the Levelling
Algorithm. This procedure can be described as follows: Starting with $%
z_{1}=z $ compute $z_{2}=z_{1}-Az_{1},$ $z_{3}=z_{2}-Bz_{2},$ $%
z_{4}=z_{3}-Az_{3},$ and so forth. Clearly, $z-z_{n}\in U+V$ and the
sequence $\{\left\Vert z_{n}\right\Vert \}_{n=1}^{\infty }$ is
nonincreasing. The question is if and when $\left\Vert z_{n}\right\Vert $
converges to the error of approximation from $U+V$?

In 1933, von Neumann \cite{19} proved that in Hilbert spaces the levelling
algorithm is always effective for any pair of closed subspaces. But in
general, for Banach spaces, this method needs additional favorable
conditions to be effective. A general result of M.Golomb \cite{6} (see also
Light and Cheney \cite[p.57]{14}) is formulated as follows

\bigskip

\textbf{Theorem 1.1 (Golomb \cite{6}).} \textit{Let $U$ and $V$ be subspaces
of a Banach space having central proximity maps $A$ and $B$, respectively,
and $U+V$ be closed. Then the sequence $\{z_{n}\}$ converges in norm to the
error of approximation from $U+V,$ that is, $\left\Vert z_{n}\right\Vert
\downarrow dist(z,U+V).$}

\bigskip

A proximity map $A$ from a Banach space $E$ onto a subspace $U$ is called a
central proximity map if for all $z\in E$ and $x\in U,$

\begin{equation*}
\left\Vert z-Az+x\right\Vert =\left\Vert z-Az-x\right\Vert \text{ (see \cite[%
Chapter 4]{14}).}
\end{equation*}

It should be noted that in a uniformly convex and uniformly smooth Banach
space it was proved by Deutsch \cite{D} that the above algorithm converges
for two closed subspaces $U$ and $V$ if the sum $U+V$ is closed. He also
noted that in any normed linear space that is not smooth one can always
construct two linear subspaces for which the algorithm does not converge.

In the Banach space setting, the levelling algorithm goes under various
names. It is fairly often called the Diliberto-Straus algorithm. Diliberto
and Straus \cite{4} were the first to consider this algorithm in the space
of continuous functions. They proved that for the problem of uniform
approximation of a bivariate function, defined on a rectangle with sides
parallel to the coordinate axes, by sums of univariate functions, the
sequence produced by the levelling algorithm converges to the desired
quantity. Later this algorithm was generalized to continuous functions
defined on a Cartesian product of two compact spaces (see \cite{13,14}).
M.Golomb \cite{6} observed that the centrality property of proximity maps
plays a key role in the analysis of Diliberto and Straus. It should be
remarked that the concept \textquotedblleft central proximity map" is due to
Golomb. We refer the readers to the monographs by Light, Cheney \cite{14}
and Khavinson \cite{12} for interesting discussions around the
Diliberto-Straus algorithm.

Let $X$ be a compact Hausdorff space. In the current paper we consider the
Diliberto-Straus algorithm in the problem of approximating from a sum of two
closed subalgebras of $C(X)$ that contain the constants. Under mild and
reasonable assumptions, we prove that the sequence produced by the levelling
algorithm converges to the error of approximation.

\bigskip

\bigskip

\begin{center}
{\large \textbf{2. Main Result}}
\end{center}

Let $X$ be a compact Hausdorff space, $C(X)$ be the space of real-valued
continuous functions on $X$ and $A_{1}\subset C(X),$ $A_{2}\subset C(X)$ be
two closed algebras that contain the constants. Define the equivalence
relation $R_{i},$ $i=1,2,$ for elements in $X$ by setting

\begin{equation*}
a\overset{R_{i}}{\sim }b\text{ if }f(a)=f(b)\text{ for all }f\in A_{i}.\eqno%
(2.1)
\end{equation*}

Then, for each $i=1,2,$ the quotient space $X_{i}=X/R_{i}$ with respect to
the relation $R_{i}$, equipped with the quotient space topology, is compact.
In addition, the natural projections $s:X\rightarrow X_{1}$ and $%
p:X\rightarrow X_{2}$ are continuous. Note that the quotient spaces $X_{1}$
and $X_{2}$ are not only compact but also Hausdorff (see, e.g., \cite[p.54]%
{12}). Obviously, in view of the Stone-Weierstrass theorem,
\begin{eqnarray*}
A_{1} &=&\{f(s(x)):~f\in C(X_{1})\}, \\
A_{2} &=&\{g(p(x)):~g\in C(X_{2})\}.
\end{eqnarray*}

In this paper, we consider the problem of approximation of a function $h\in
C(X)$ from the set $A_{1}+A_{2}.$ More precisely, we are interested in
algorithmic methods for computing the error of approximation

\begin{equation*}
E(h)=\inf_{w\in A_{1}+A_{2}}\left\Vert h-w\right\Vert .
\end{equation*}

In the sequel, we assume that the algebras $A_{1}$ and $A_{2}$ obey the
following property, which we call the $C$-property: For any function $h\in
C(X),$ the real functions

\begin{eqnarray*}
f_{1}(a) &=&\max_{\substack{ x\in X  \\ s(x)=a}}h(x),\text{ }f_{2}(a)=\min
_{\substack{ x\in X  \\ s(x)=a}}h(x),\text{ }a\in X_{1}, \\
g_{1}(b) &=&\max_{\substack{ x\in X  \\ p(x)=b}}h(x),\text{ }g_{2}(b)=\min
_{\substack{ x\in X  \\ p(x)=b}}h(x),\text{ }b\in X_{2}
\end{eqnarray*}%
are continuous on the appropriate sets $X_{1}$ and $X_{2}.$ Note that for a
given compact space $X$ many subalgebras of $C(X)$ do not possess the $C$%
-property. For example, in the above special case of $\mathbb{R}^{2},$ the
above $\max $ functions are not continuous if $X=[0,1]\times \lbrack
0,2]\cup \lbrack 1,2]\times \lbrack 0,1]$ and $h(x,y)=xy.$ In order for the $%
C$-property to be fulfilled, the natural quotient mappings $s(x)$ and $p(x)$
should satisfy certain conditions. The following proposition provides a
sufficient condition of such type for compact sequential spaces. A
sequential space is a topological space with the property that a set is open
iff every sequence $x_{n}$ converging to a point in the set is, itself,
eventually in the set (i.e. there exists $N$ such that $x_{n}$ is in the set
for all $n\geq N$). Note that some essential properties and full
characterization of sequential spaces via quotient mappings were given in
the fundamental papers of Franklin (see \cite{Fr1,Fr2}).

\bigskip

\textbf{Proposition 2.1.} \textit{Let $X$ be a compact sequential Hausdorff
space and $A$ be a closed subalgebra of $C(X)$ that contains the constants.
Let $X_{1}$ be a quotient space generated by the equivalence relation (2.1)
and $s:X\rightarrow X_{1}$ be the natural quotient mapping. Then the
functions $f_{1}$ and $f_{2}$ are continuous on $X_{1}$ for any $h\in C(X)$
if for any two points $x$ and $y$ with $s(x)=s(y)$ and any sequence $%
\{x_{n}\}_{n=1}^{\infty }$ tending to $x,$ there exists a sequence $%
\{y_{n}\}_{n=1}^{\infty }$ tending to $y$ such that $s(y_{n})=s(x_{n}),$ for
all $n=1,2,...$}

\begin{proof} Suppose the contrary. Suppose that the above hypothesis on the
quotient mapping $s$ holds, but one of the functions $f_{1}$ and $f_{2}$ is
not continuous. Without loss of generality assume that $f_{1}$ is not
continuous on $X_{1}$. Let $f_{1}$ be discontinuous at a point $a_{0}\in
X_{1}.$ Note that by the result of Franklin \cite[Proposition 1.2]{Fr1}, a
quotient image of a sequential space is sequential. Thus the quotient space $%
X_{1}$ is sequential. Then there exists a number $\varepsilon >0$ and a
sequence $\{a_{n}\}_{n=1}^{\infty }\subset X_{1}$ tending to $a_{0},$ such
that

\begin{equation*}
\left\vert f_{1}(a_{n})-f_{1}(a_{0})\right\vert >\varepsilon ,\eqno(2.2)
\end{equation*}%
for all $n=1,2,...$. Since the function $h$ is continuous on $X$, there
exist points $x_{n}\in X$, $n=0,1,2,...,$ such that $h(x_{n})=f_{1}(a_{n}),$
$s(x_{n})=a_{n},$ for $n=0,1,2,...$. \ Thus the inequality (2.2) can be
written as

\begin{equation*}
\left\vert h(x_{n})-h(x_{0})\right\vert >\varepsilon ,\eqno(2.3)
\end{equation*}%
for all $n=1,2,...$. Since $X$ is compact and sequential, it is sequentially
compact (see \cite[Theorem 3.10.31]{En}); hence the sequence $%
\{x_{n}\}_{n=1}^{\infty }$ has a convergent subsequence. Without loss of
generality assume that $\{x_{n}\}_{n=1}^{\infty }$ itself converges to a
point $y_{0}\in X.$ Then $s(x_{n})\rightarrow s(y_{0}),$ as $n\rightarrow
\infty .$ But by the assumption, we also have $s(x_{n})\rightarrow s(x_{0}),$
as $n\rightarrow \infty .$ Therefore, since $X_{1}$ is Hausdorff, $%
s(y_{0})=s(x_{0})=a_{0}.$ Note that $x_{0}$ and $y_{0}$ cannot be the same
point, since the equality $x_{0}=y_{0}$ violates the condition (2.3). By the
hypothesis of the proposition, we must have a sequence $\{z_{n}\}_{n=1}^{%
\infty }$ such that $z_{n}\rightarrow x_{0}$ and

\begin{equation*}
s(z_{n})=s(x_{n}),
\end{equation*}%
for all $n=1,2,...$. Since $s(x_{n})=a_{n},$ $n=1,2,...$, and on each level
set $\{x\in X:s(x)=a_{n}\},$ the function $h$ takes its maximum value at $%
x_{n}$ we obtain that

\begin{equation*}
h(z_{n})\leq h(x_{n}),\text{ }n=1,2,...
\end{equation*}%
Taking the limit in the last inequality as $n\rightarrow \infty ,$ gives us
the new inequality

\begin{equation*}
h(x_{0})\leq h(y_{0}).\eqno(2.4)
\end{equation*}

Recall that on the level set $\{x\in X:$ $s(x)=a_{0}\}$, the function $h$
takes its maximum at $x_{0}$. Thus from (2.4) we conclude that $%
h(x_{0})=h(y_{0}).$ This last equality contradicts the choice of the
positive $\varepsilon $ in (2.3), since $h(x_{n})\rightarrow h(y_{0}),$ as $%
n\rightarrow \infty .$ The obtained contradiction shows that the function $%
f_{1}$ is continuous on $X_{1}.$ In the same way one can prove that $f_{2}$
is continuous on $X_{1}.$ \end{proof}

\bigskip

\textbf{Example 1.} It is not difficult to see that the inner product
function $s(\mathbf{x})=\mathbf{a}\cdot \mathbf{x}$, where $\mathbf{a}$ is a
nonzero vector in $\mathbb{R}^{d}$, satisfies the hypothesis of Proposition
2.1, if $\mathbf{x}$ varies in a compact convex set $Q\subset \mathbb{R}^{d}$%
. Indeed, let $\mathbf{x}_{0}$ and $\mathbf{y}_{0}$ be any two points in $Q$
such that $s(\mathbf{x}_{0})=s(\mathbf{y}_{0}).$ Take any sequence $\{%
\mathbf{x}_{n}\}_{n=1}^{\infty }\subset Q$, $\mathbf{x}_{n}\rightarrow
\mathbf{x}_{0},$ as $n\rightarrow \infty .$ We must show that there exists a
sequence $\{\mathbf{y}_{n}\}_{n=1}^{\infty }\subset Q$ with the properties
that $\mathbf{y}_{n}\rightarrow \mathbf{y}_{0},$ as $n\rightarrow \infty $,
and $s(\mathbf{x}_{n})=s(\mathbf{y}_{n}),$ for $n=1,2,...$ To show this,
first note that we have points $\mathbf{z}_{1},\mathbf{z}_{2}\in Q$ such
that $s(\mathbf{x}_{n})\in \lbrack s(\mathbf{z}_{1}),s(\mathbf{z}_{2})],$
for all $n=1,2,...$ One of, or both, the numbers $s(\mathbf{z}_{1})$ and $s(%
\mathbf{z}_{2})$ may equal to $s(\mathbf{x}_{0})$, and in this case we
assume that $\mathbf{z}_{1}$ and/or $\mathbf{z}_{2}$ coincides with $\mathbf{%
y}_{0}$. Consider now the line segments $[\mathbf{z}_{1},\mathbf{y}_{0}]$
and $[\mathbf{y}_{0},\mathbf{z}_{2}]$ (which may be degenerated into the
point $\mathbf{y}_{0}$). Set $L=[\mathbf{z}_{1},\mathbf{y}_{0}]\cup \lbrack
\mathbf{y}_{0},\mathbf{z}_{2}].$ Since $Q$ is convex, $L\subset Q.$ The
function $s$ is continuous on $L.$ Hence by the Intermediate Value Theorem,
there exists a sequence $\{\mathbf{y}_{n}\}_{n=1}^{\infty }\subset L$ such
that $s(\mathbf{y}_{n})=s(\mathbf{x}_{n}),$ for $n=1,2,...$ Since $s(\mathbf{%
x}_{n})\rightarrow s(\mathbf{x}_{0})$ and $s(\mathbf{x}_{0})=s(\mathbf{y}%
_{0})$, it follows that $s(\mathbf{y}_{n})\rightarrow s(\mathbf{y}_{0}),$ as
$n\rightarrow \infty $. Now from this fact and the inclusion $\{\mathbf{y}%
_{n}\}_{n=1}^{\infty }\subset L$ we easily derive that $\mathbf{y}%
_{n}\rightarrow \mathbf{y}_{0},$ as $n\rightarrow \infty $.

\bigskip

\textbf{Example 2.} Let $Q$ be a compact convex set in $\mathbb{R}^{2}$.
From Example 1 it follows that the $C$-property holds for the algebras of
univariate functions $A_{1}=\{f(x):$ $f\in Q_{x}\}$ and $A_{2}=\{g(y):$ $%
g\in Q_{y}\},$ where $Q_{x}$ and $Q_{y}$ are projections of $Q$ into the
coordinate axes $x$ and $y$, respectively.

\bigskip

\textbf{Example 3.} The hypothesis of Proposition 2.1 strictly depends on
the considered space $X.$ That is, the natural quotient mapping $s$ may
satisfy this hypothesis, but for many closed subsets $E\subset X,$ it may
happen that the restriction of $s$ to $E$ no longer satisfies it. For
example, let $K$ be the unit square in the $xy$ plane and $K_{1}=[0,1]\times
\lbrack 0,\frac{1}{2}]\cup \lbrack 0,\frac{1}{2}]\times \lbrack 0,1].$
Consider the algebra $U$ of univariate functions depending only on the
variable $x.$ Clearly, the coordinate projection $s(x,y)=x$ satisfies the
hypothesis if $s$ is considered over $K$. This is not true if $s$ is
considered over the set $K_{1}.$ Indeed, for the sequence $\{(\frac{1}{2}+%
\frac{1}{n+1},\frac{1}{2})\}_{n=1}^{\infty }\subset K_{1},$ which tends to $(%
\frac{1}{2},\frac{1}{2}),$ we cannot find a sequence $\{(x_{n},y_{n})%
\}_{n=1}^{\infty }\subset K_{1}$ tending to $(\frac{1}{2},1)$ such that $%
x_{n}=\frac{1}{2}+\frac{1}{n+1},$ $n=1,2,...$

\bigskip

Define the following operators
\begin{equation*}
F:C(X)\rightarrow A_{1},~~Fh(a)=\frac{1}{2}\left( \max_{\substack{ x\in X
\\ s(x)=a}}h(x)+\min_{\substack{ x\in X  \\ s(x)=a}}h(x)\right) ,\text{ \ \
for all }a\in X_{1},
\end{equation*}%
and%
\begin{equation*}
G:C(X)\rightarrow A_{2},~~Gh(b)=\frac{1}{2}\left( \max_{\substack{ x\in X
\\ p(x)=b}}h(x)+\min_{\substack{ x\in X  \\ p(x)=b}}h(x)\right) ,\text{ \ \
for all }b\in X_{2}.
\end{equation*}

Since the algebras $A_{i},$ $i=1,2,$ both enjoy the $C$-property, for each
function $h\in C(X)$, the functions $Fh$ and $Gh$ are continuous on $X_{1}$
and $X_{2}$, respectively. Via the quotient mappings $s$ and $p,$ these
functions $Fh$ and $Gh$ can be considered also as functions defined on $X$.
Since $s$ and $p$ are continuous on $X$, the functions $Fh$ and $Gh$ are
continuous on $X$ and hence belong to the algebras $A_{1}$ and $A_{2}$
respectively.

The following theorem plays a key role in the proof of our main result
(Theorem 2.3).

\bigskip

\textbf{Theorem 2.2.} \textit{Assume $X$ is a compact Hausdorff space and $%
A_{i},$ $i=1,2,$ are closed subalgebras of $C(X)$ that contain the
constants. In addition, assume that the $C$-property holds for these
subalgebras. Then the operators $F$ and $G$ are central proximity maps onto $%
A_{1}$ and $A_{2}$, respectively. In addition, these operators are
non-expansive. That is,}
\begin{equation*}
\left\Vert Fv_{1}-Fv_{2}\right\Vert \leq \left\Vert v_{1}-v_{2}\right\Vert
\text{ \textit{and} }\left\Vert Gv_{1}-Gv_{2}\right\Vert \leq \left\Vert
v_{1}-v_{2}\right\Vert ,
\end{equation*}%
\textit{for all $v_{1},v_{2}\in C(X)$.}

\begin{proof} We prove this theorem for the operator $F.$ A proof for $G$ can
be carried out in the same way.

Clearly, on the level set $\{x\in X:$ $s(x)=a\},$ the constant $(Fh)(a)$ is
a best approximation to $h,$ among all constants. Varying over $a\in X_{1},$
we obtain a best approximating function $Fh:X\rightarrow \mathbb{R}$, which
is, due to the $C$-property, in the algebra $A_{1}.$

Now let us prove that the proximity map $F$ is a central proximity map. In
other words, we must prove that for any functions $h\in C(X)$ and $%
f=f(s(x))\in A_{1}$,

\begin{equation*}
\left\Vert h-Fh-f\right\Vert =\left\Vert h-Fh+f\right\Vert .\eqno(2.5)
\end{equation*}

Put $u=h-Fh.$ There exists a point $x_{0}\in X$ such that

\begin{equation*}
\left\Vert u+f\right\Vert =\left\vert u(x_{0})+f(s(x_{0}))\right\vert .
\end{equation*}

First assume that $\left\vert u(x_{0})+f(s(x_{0}))\right\vert
=u(x_{0})+f(s(x_{0})).$ Since $Fu=0,$

\begin{equation*}
\max_{\substack{ x\in X  \\ s(x)=a}}u(x)=-\min_{\substack{ x\in X  \\ s(x)=a
}}u(x),~\text{for all }a\in X_{1}.\eqno(2.6)
\end{equation*}

Let

\begin{equation*}
\min_{\substack{ x\in X  \\ s(x)=s(x_{0})}}u(x)=u(x_{1}).\eqno(2.7)
\end{equation*}%
Then from (2.6) and (2.7) it follows that

\begin{equation*}
-u(x_{1})\geq u(x_{0}).
\end{equation*}%
Taking the last inequality and the equality $s(x_{1})=s(x_{0})$ into account
we may write

\begin{equation*}
\left\Vert u-f\right\Vert \geq f(s(x_{1}))-u(x_{1})\geq
f(s(x_{0}))+u(x_{0})=\left\Vert u+f\right\Vert .\eqno(2.8)
\end{equation*}%
Changing in (2.8) the function $f$ to $-f$ gives the reverse inequality $%
\left\Vert u+f\right\Vert \geq \left\Vert u-f\right\Vert .$ Thus (2.5) holds.

Note that if $\left\vert u(x_{0})+f(s(x_{0}))\right\vert
=-(u(x_{0})+f(s(x_{0})))$, then by replacing Eq (2.7) by

\begin{equation*}
\max_{\substack{ x\in X  \\ s(x)=s(x_{0})}}u(x)=u(x_{1}).\eqno(2.9)
\end{equation*}%
we will derive from (2.6) and (2.9) that $\ u(x_{1})\geq -u(x_{0}).$ This
inequality is then used to obtain the estimation

\begin{equation*}
\left\Vert u-f\right\Vert \geq -(f(s(x_{1})-u(x_{1}))\geq
-(f(s(x_{0}))+u(x_{0}))=\left\Vert u+f\right\Vert ,
\end{equation*}%
which in turn yields (2.5). The centrality has been proven.

Now we prove that the operator $F$ is non-expansive. First note that it is
nondecreasing. That is, if for all $x\in X,$ $h_{1}(x)\leq h_{2}(x),$ then $%
Fh_{1}(s(x))\leq Fh_{2}(s(x))$. Besides, $F(h+c)=Fh+c$, for any real number $%
c$. Let now $v_{1}$ and $v_{2}$ be arbitrary functions in $C(X)$. Put $%
c=\left\Vert v_{1}-v_{2}\right\Vert .$ Then for any $x\in X$, we can write

\begin{equation*}
v_{2}(x)-c\leq v_{1}(x)\leq v_{2}(x)+c
\end{equation*}%
and further

\begin{equation*}
Fv_{2}(s(x))-c\leq Fv_{1}(s(x))\leq Fv_{2}(s(x))+c.
\end{equation*}%
From the last inequality we obtain that

\begin{equation*}
\left\Vert Fv_{1}-Fv_{2}\right\Vert \leq c=\left\Vert v_{1}-v_{2}\right\Vert
.
\end{equation*}%
Thus the operator $F$ is non-expansive. \end{proof}

\smallskip

Consider the iterations

\begin{equation*}
h_{1}(x)=h(x),\text{ }h_{2n}=h_{2n-1}-Fh_{2n-1},\text{ }%
h_{2n+1}=h_{2n}-Gh_{2n},\text{ }n=1,2,...\text{.}
\end{equation*}%
From Theorem 2.2 and the above-mentioned general result of Golomb (see
Theorem 1.1) one can obtain the following theorem.

\bigskip

\textbf{Theorem 2.3.} \textit{Let all the assumptions of Theorem 2.2 hold
and $A_{1}+A_{2}$ be closed in $C(X).$ Then $\left\Vert h_{n}\right\Vert $
converges to the error of approximation $E(h).$}

\bigskip

Does the Diliberto and Straus algorithm converge without the closedness
assumption on the sum $A_{1}+A_{2}$? We do not yet know a complete answer to
this question. Recall that we obtain Theorem 2.3 by using the general result
of Golomb (see Theorem 1.1), in which the \textquotedblleft closedness" is a
major hypothesis. Below we give a different (and independent of Golomb's
result) proof of this theorem, where it is shown how we use the mentioned
closedness. It should be remarked that in the simplest case of approximation
by sums of univariate functions, the known classical proofs of the
Diliberto-Straus algorithm require that one could form a closed bolt (for
this terminology see below) from an arbitrary bolt $(x_{1},...,x_{n})$ of a
given compact set $Q\subset \mathbb{R}^{2}$ by adding a point $y\in Q$,
whose first and second coordinates are equal to that of $x_{1}$ and $x_{n}$
respectively (see, e.g., \cite{4,12,13}). Note that this point may not lie
in $Q,$ unless $Q$ is a Cartesian product of two compact sets. The main idea
behind our proof is to make use of weak$^{\text{*}}$ cluster points of some
sequence of \textquotedblleft unclosed bolt functionals" instead of
considering only \textquotedblleft closed bolt functionals". We hope that
this idea can be useful in future attempts to prove the convergence of the
Diliberto-Straus algorithm without the closedness assumption.

\bigskip

For the further analysis we need the following objects called
\textquotedblleft bolts of lightning" or simply \textquotedblleft bolts".

\bigskip

\textbf{Definition 2.2 (see \cite{Mar1,Mar2}).} \textit{Let $X$ be a compact
space and $A_{i}$, $i=1,2,$ be subalgebras of $C(X)$ that contain the
constants. Let $s$ and $p$ be the natural quotient mappings generated by the
equivalence relation (2.1). A finite ordered set $l=\{x_{1},x_{2},...,x_{n}%
\}\subset X$, where $x_{i}\neq x_{i+1}$, with either $%
s(x_{1})=s(x_{2}),p(x_{2})=p(x_{3}),s(x_{3})=s(x_{4}),...$ or $%
p(x_{1})=p(x_{2}),s(x_{2})=s(x_{3}),p(x_{3})=p(x_{4}),...$ is called a bolt
with respect to $(A_{1},A_{2})$.}

\bigskip

If in a bolt $\{x_{1},...,x_{n},x_{n+1}\}$, $x_{n+1}=x_{1}$ and $n$ is an
even number, then the bolt $\{x_{1},...,x_{n}\}$ is said to be closed.
Bolts, in the case when $X\subset \mathbb{R}^{2}$ and the algebras $A_{1}$
and $A_{2}$ coincide with the spaces of univariate functions $\varphi (x)$
and $\psi (y),$ respectively, are geometrically explicit objects. In this
case, a bolt is a finite ordered set $\{x_{1},x_{2},...,x_{n}\}$ in $\mathbb{%
R}^{2}$ with the line segments $[x_{i},x_{i+1}],$ $i=1,...,n,$ alternatively
perpendicular to the $x$ and $y$ axes (see, e.g., \cite%
{Arn,Is1,Is2,Is3,12,Mar1}). These objects were first introduced by Diliberto
and Straus \cite{4} (in \cite{4}, they are called \textquotedblleft
permissible lines\textquotedblright ). They appeared further in a number of
papers with several different names such as \textquotedblleft paths" (see,
e.g., \cite{13,14}), \textquotedblleft trips\textquotedblright\ (see \cite%
{Mar2,Med}), \textquotedblleft links" (see, e.g., \cite{Cow,Klo1,Klo2}),
etc. The term \textquotedblleft bolt of lightning" was due to Arnold \cite%
{Arn}.

With each bolt $l=\{x_{1},...,x_{n}\}$ with respect to $(A_{1},A_{2})$, we
associate the following bolt functional

\begin{equation*}
r_{l}(h)=\frac{1}{n}\sum_{i=1}^{n}(-1)^{n+1}h(x_{i}).
\end{equation*}

It is an exercise to check that $r_{l}$ is a linear bounded functional on $%
C(X)$ with the norm $\left\Vert r_{l}\right\Vert \leq 1$ and $\left\Vert
r_{l}\right\Vert =1$ if and only if the set of points $x_{i}$ with odd
indices $i$ does not intersect with the set of points with even indices.
Besides, if $l$ is closed, then $r_{l}\in (A_{1}+A_{2})^{\perp },$ where $%
(A_{1}+A_{2})^{\perp }$ is the annihilator of the subspace $%
A_{1}+A_{2}\subset C(X).$ If $l$ is not closed, then $r_{l}$ is generally
not an annihilating functional. However, it obeys the following important
inequality

\begin{equation*}
\left\vert r_{l}(f)\right\vert \leq \frac{2}{n}\left\Vert f\right\Vert ,\eqno%
(2.10)
\end{equation*}%
for all $f\in A_{i}$, $i=1,2$. This inequality means that for bolts $l$ with
sufficiently large number of points, $r_{l}$ behaves like an annihilating
functional.

Now we are able to prove Theorem 2.3.

\begin{proof} First, let us write the above iteration in the following form
\begin{eqnarray*}
h_{1} &=&h,\text{ }h_{n+1}=h_{n}-q_{n},\text{ where} \\
q_{n} &=&Fh_{n},\text{ if }n\text{ is odd;} \\
q_{n} &=&Gh_{n},\text{ if }n\text{ is even.}
\end{eqnarray*}

Introduce the functions
\begin{eqnarray*}
u_{n} &=&q_{1}+\cdot \cdot \cdot +q_{2n-1}, \\
v_{0} &=&0,\text{ }v_{n}=q_{2}+\cdot \cdot \cdot +q_{2n},\text{ }n=1,2,...
\end{eqnarray*}%
Clearly, $u_{n}\in A_{1}$ and $v_{n}\in A_{2}.$ Besides, $%
h_{2n}=h-u_{n}-v_{n-1}$ and $h_{2n+1}=h-u_{n}-v_{n},$ for $n=1,2,...$.

It is easy to see that the following inequalities hold

\begin{equation*}
\left\Vert h_{1}\right\Vert \geq \left\Vert h_{2}\right\Vert \geq \left\Vert
h_{3}\right\Vert \geq \cdot \cdot \cdot \geq E(h).
\end{equation*}%
Therefore, there exists the limit%
\begin{equation*}
M=\lim_{n\rightarrow \infty }\left\Vert h_{n}\right\Vert \geq E(h).
\end{equation*}

It is a consequence of the Hahn-Banach extension theorem that

\begin{equation*}
E(h)=\sup_{\substack{ r\in (A_{1}+A_{2})^{\bot } \\ \left\Vert r\right\Vert
\leq 1}}\left\vert r(h)\right\vert ,
\end{equation*}%
where the $\sup $ is attained by some functional. To complete the proof it
is enough to show that for any $\varepsilon >0,$ there exists a functional $%
r_{0}$ such that $r_{0}\in (A_{1}+A_{2})^{\bot }$, $\left\Vert
r_{0}\right\Vert \leq 1$ and $\left\vert r_{0}(h)\right\vert \geq
M-\varepsilon .$

Let $\varepsilon $ be an arbitrarily small positive real number. For each
positive integer $k=1,2,...,$ set $\delta _{k}=\frac{\varepsilon }{2^{2k}}.$
There exists a number $n_{k}$ such that for all $n\geq n_{k}$,%
\begin{equation*}
\left\Vert h_{n}\right\Vert \leq M+\delta _{k}.
\end{equation*}

Without loss of generality we may assume that $n_{k}$ is even. In the
following, for each $k$, we are going to construct a bolt $%
l_{k}=\{x_{1}^{k},...,x_{2k+1}^{k}\}$ with the property that $\left\vert
r_{l_{k}}(h_{n_{k}})\right\vert \geq M-\varepsilon $. This will lead us to
the above mentioned functional $r_{0}$.

Since $Fh_{2m}=0,$ for $m=1,2,...,$ and $n_{k}+2k$ is an even number,

\begin{equation*}
\max_{\substack{ x\in X \\ s(x)=a}}h_{n_{k}+2k}(x)=-\min_{\substack{ x\in X
\\ s(x)=a}}h_{n_{k}+2k}(x),\text{ \ \ for all }a\in X_{1}.
\end{equation*}%
Then there exist points $x_{1}$ and $x_{2}$ such that $h_{n_{k}+2k}(x_{1})=%
\left\Vert h_{n_{k}+2k}\right\Vert $, $h_{n_{k}+2k}(x_{2})=-\left\Vert
h_{n_{k}+2k}\right\Vert $ and $s(x_{1})=s(x_{2}).$ This can be written in
the form

\begin{equation*}
h_{n_{k}+2k-1}(x_{1})-Fh_{n_{k}+2k-1}(s(x_{1}))=\left\Vert
h_{n_{k}+2k}\right\Vert ;\eqno(2.11)
\end{equation*}

\begin{equation*}
h_{n_{k}+2k-1}(x_{2})-Fh_{n_{k}+2k-1}(s(x_{2}))=-\left\Vert
h_{n_{k}+2k}\right\Vert .\eqno(2.12)
\end{equation*}%
Therefore,

\begin{equation*}
Fh_{n_{k}+2k-1}(s(x_{1}))=Fh_{n_{k}+2k-1}(s(x_{2}))=
\end{equation*}%
\begin{equation*}
=h_{n_{k}+2k-1}(x_{1})-\left\Vert h_{n_{k}+2k}\right\Vert \leq \left\Vert
h_{n_{k}+2k-1}\right\Vert -M\leq \delta _{k}
\end{equation*}%
and

\begin{equation*}
Fh_{n_{k}+2k-1}(s(x_{1}))=Fh_{n_{k}+2k-1}(s(x_{2}))=
\end{equation*}%
\begin{equation*}
=h_{n_{k}+2k-1}(x_{2})+\left\Vert h_{n_{k}+2k}\right\Vert \geq -(\left\Vert
h_{n_{k}+2k-1}\right\Vert -M)\geq -\delta _{k}.
\end{equation*}%
That is,

\begin{equation*}
-\delta _{k}\leq Fh_{n_{k}+2k-1}(s(x_{1}))=Fh_{n_{k}+2k-1}(s(x_{2}))\leq
\delta _{k}\eqno(2.13)
\end{equation*}

From (2.11)-(2.13) we obtain that

\begin{equation*}
h_{n_{k}+2k-1}(x_{1})=Fh_{n_{k}+2k-1}(s(x_{1}))+\left\Vert
h_{n_{k}+2k}\right\Vert \geq M-\delta _{k};\eqno(2.14)
\end{equation*}

\begin{equation*}
h_{n_{k}+2k-1}(x_{2})=Fh_{n_{k}+2k-1}(s(x_{1}))-\left\Vert
h_{n_{k}+2k}\right\Vert \leq -M+\delta _{k}.\eqno(2.15)
\end{equation*}%
Since $Gh_{2m+1}=0,$ for $m=1,2,...,$

\begin{equation*}
\max_{\substack{ x\in X  \\ p(x)=b}}h_{n_{k}+2k-1}(x)=-\min_{\substack{ x\in
X  \\ p(x)=b}}h_{n_{k}+2k-1}(x),\text{ \ \ for all }b\in X_{2}.
\end{equation*}%
Then from (2.14) and (2.15) it follows that there exist points $x_{3}$ and $%
x_{3}^{^{\prime }}$ satisfying the relations

\begin{equation*}
h_{n_{k}+2k-1}(x_{3})\geq M-\delta _{k},\text{ \ }p(x_{3})=p(x_{2})\eqno%
(2.16)
\end{equation*}%
and

\begin{equation*}
h_{n_{k}+2k-1}(x_{3}^{^{\prime }})\leq -M+\delta _{k},\text{ \ }%
p(x_{3}^{^{\prime }})=p(x_{1}).\eqno(2.17)
\end{equation*}

In the inequalities (2.14)-(2.17), replace $h_{n_{k}+2k-1}$ by $%
h_{n_{k}+2k-2}-Gh_{n_{k}+2k-2}.$ Then we have the following estimates.

\begin{equation*}
h_{n_{k}+2k-2}(x_{1})-Gh_{n_{k}+2k-2}(p(x_{1}))\geq M-\delta _{k};\eqno(2.18)
\end{equation*}

\begin{equation*}
h_{n_{k}+2k-2}(x_{2})-Gh_{n_{k}+2k-2}(p(x_{2}))\leq -M+\delta _{k};\eqno%
(2.19)
\end{equation*}

\begin{equation*}
h_{n_{k}+2k-2}(x_{3})-Gh_{n_{k}+2k-2}(p(x_{3}))\geq M-\delta _{k};\eqno(2.20)
\end{equation*}

\begin{equation*}
h_{n_{k}+2k-2}(x_{3}^{\prime })-Gh_{n_{k}+2k-2}(p(x_{3}^{^{\prime }}))\leq
-M+\delta _{k}.\eqno(2.21)
\end{equation*}%
The following inequalities are obvious.
\begin{eqnarray*}
h_{n_{k}+2k-2}(x_{1}) &\leq &\left\Vert h_{n_{k}+2k-2}\right\Vert \leq
M+\delta _{k}, \\
h_{n_{k}+2k-2}(x_{3}) &\leq &\left\Vert h_{n_{k}+2k-2}\right\Vert \leq
M+\delta _{k}, \\
h_{n_{k}+2k-2}(x_{2}) &\geq &-\left\Vert h_{n_{k}+2k-2}\right\Vert \geq
-M-\delta _{k}, \\
h_{n_{k}+2k-2}(x_{3}^{^{\prime }}) &\geq &-\left\Vert
h_{n_{k}+2k-2}\right\Vert \geq -M-\delta _{k}.
\end{eqnarray*}%
Considering these obvious inequalities in (2.18)-(2.21), we can write that

\begin{equation*}
-2\delta _{k}\leq
Gh_{n_{k}+2k-2}(p(x_{1}))=Gh_{n_{k}+2k-2}(p(x_{3}^{^{\prime }}))\leq 2\delta
_{k};\eqno(2.22)
\end{equation*}

\begin{equation*}
-2\delta _{k}\leq Gh_{n_{k}+2k-2}(p(x_{2}))=Gh_{n_{k}+2k-2}(p(x_{3}))\leq
2\delta _{k}.\eqno(2.23)
\end{equation*}%
Taking into account (2.22) and (2.23) in the estimates (2.18)-(2.21), we
obtain that
\begin{eqnarray*}
h_{n_{k}+2k-2}(x_{1}) &\geq &M-3\delta _{k}, \\
h_{n_{k}+2k-2}(x_{2}) &\leq &-M+3\delta _{k},
\end{eqnarray*}

\begin{equation*}
h_{n_{k}+2k-2}(x_{3})\geq M-3\delta _{k},\eqno(2.24)
\end{equation*}

\begin{equation*}
h_{n_{k}+2k-2}(x_{3}^{^{\prime }})\leq -M+3\delta _{k}.\eqno(2.25)
\end{equation*}

Now since

\begin{equation*}
\max_{\substack{ x\in X  \\ s(x)=a}}h_{n_{k}+2k-2}(x)=-\min_{\substack{ x\in
X  \\ s(x)=a}}h_{n_{k}+2k-2}(x),\text{ \ \ for all }a\in X_{1},
\end{equation*}%
from (2.24) and (2.25) it follows that there exist points $x_{4}$ and $%
x_{4}^{^{\prime }}$ satisfying

\begin{equation*}
h_{n_{k}+2k-2}(x_{4})\leq -M+3\delta _{k},\text{ }s(x_{4})=s(x_{3}),
\end{equation*}%
and

\begin{equation*}
h_{n_{k}+2k-2}(x_{4}^{^{\prime }})\geq M-3\delta _{k},\text{ }%
s(x_{4}^{^{\prime }})=s(x_{2}).
\end{equation*}

Repeating the above process for the function $%
h_{n_{k}+2k-2}(x)=h_{n_{k}+2k-3}(x)-Fh_{n_{k}+2k-3}(s(x))$ we obtain that

\begin{equation*}
-4\delta _{k}\leq Fh_{n_{k}+2k-3}(s(x_{i}))\leq 4\delta _{k},\text{ }%
i=1,2,3,4,\eqno(2.26)
\end{equation*}%
and
\begin{eqnarray*}
h_{n_{k}+2k-3}(x_{1}) &\geq &M-7\delta _{k}; \\
h_{n_{k}+2k-3}(x_{2}) &\leq &-M+7\delta _{k}; \\
h_{n_{k}+2k-3}(x_{3}) &\geq &M-7\delta _{k}; \\
h_{n_{k}+2k-3}(x_{4}) &\leq &-M+7\delta _{k}; \\
h_{n_{k}+2k-3}(x_{4}^{^{\prime }}) &\geq &M-7\delta _{k}.
\end{eqnarray*}

By the same way as above, we can find a point $x_{5}$ such that $%
p(x_{5})=p(x_{4})$ and

\begin{equation*}
h_{n_{k}+2k-3}(x_{5})\geq M-7\delta _{k}.
\end{equation*}

Continuing this process until we reach the function $h_{n_{k}}$, we obtain
the points $x_{1},x_{2}$,..., $x_{2k+2}$ satisfying $s(x_{1})=s(x_{2}),$ $%
p(x_{2})=p(x_{3}),...,$ $s(x_{2k+1})=s(x_{2k+2}).$ Clearly, these points, in
the given order, form a bolt, which we denote by $l_{k}.$ Note that in the
considered process, we also deal with the points $x_{3}^{^{\prime
}},x_{4}^{^{\prime }},$ etc., but these points play only an auxiliary role:
they are needed in obtaining the inequalities (2.22), (2.23), (2.26), etc.
At the points of the bolt $l_{k}$, the values of $h_{n_{k}}$ obey the
estimates

\begin{equation*}
h_{n_{k}}(x_{i})\geq M-(2^{2k}-1)\delta _{k}\geq M-\varepsilon ,\text{ for }%
i=1,3,...,2k+1
\end{equation*}%
and

\begin{equation*}
h_{n_{k}}(x_{j})\leq -M+(2^{2k}-1)\delta _{k}\leq -M+\varepsilon ,\text{ for
}j=2,4,...,2k+2.
\end{equation*}%
Using these inequalities, we can estimate the absolute value of $%
r_{l_{k}}(h_{n_{k}})$ as follows

\begin{equation*}
\left\vert r_{l_{k}}(h_{n_{k}})\right\vert \geq \frac{(k+1)(M-\varepsilon
)-(k+1)(-M+\varepsilon )}{2k+2}=M-\varepsilon .\eqno(2.27)
\end{equation*}%
We can also write the following obvious inequality
\begin{equation*}
\left\vert r_{l_{k}}(h)\right\vert =\left\vert
r_{l_{k}}(h_{n_{k}}+u_{n_{k}/2}+v_{n_{k}/2-1})\right\vert \geq \left\vert
r_{l_{k}}(h_{n_{k}})\right\vert -\left\vert
r_{l_{k}}(u_{n_{k}/2})\right\vert -\left\vert
r_{l_{k}}(v_{n_{k}/2-1})\right\vert .\eqno(2.28)
\end{equation*}

Considering (2.27) and (2.10) in (2.28) we obtain that

\begin{equation*}
\left\vert r_{l_{k}}(h)\right\vert \geq M-\varepsilon -\frac{2}{2k+2}\left(
\left\Vert u_{n_{k}/2}\right\Vert +\left\Vert v_{n_{k}/2-1}\right\Vert
\right) .\eqno(2.29)
\end{equation*}%
\qquad \qquad\ \

At this stage we use our assumption of closedness of the subspace $%
A_{1}+A_{2}.$ Note that a subspace $B=B_{1}+B_{2}$ of a Banach space is
closed if and only if there exists a $K<\infty $ such that each element $b$
in $B$ has a representation $b=b_{1}+b_{2},$ where $b_{1}\in B_{1},$ $%
b_{2}\in B_{2}$ and $\max (\left\Vert b_{1}\right\Vert ,\left\Vert
b_{2}\right\Vert )\leq K\left\Vert b\right\Vert $ (see \cite{Mar2}). Now,
since $A_{1}+A_{2}$ is closed, there exists a constant $K$ such that
\begin{equation*}
\left\Vert u_{n}\right\Vert +\left\Vert v_{n-1}\right\Vert \leq K\left\Vert
u_{n}+v_{n-1}\right\Vert ,\text{ }n=1,2,...\eqno(2.30)
\end{equation*}%
Note that the sequence $\{\left\Vert u_{n}+v_{n-1}\right\Vert
\}_{n=1}^{\infty }$ is uniformly bounded. Indeed,

\begin{equation*}
\left\Vert u_{n}+v_{n-1}\right\Vert =\left\Vert h-h_{2n}\right\Vert \leq
2\left\Vert h\right\Vert ,\text{ for all }n=1,2,...\eqno(2.31)
\end{equation*}%
From (2.29), (2.30) and (2.31) it follows that there exists a constant $C>0$
such that%
\begin{equation*}
\left\vert r_{l_{k}}(h)\right\vert \geq M-\varepsilon -\frac{C}{2k+2}.\eqno%
(2.32)
\end{equation*}

Thus for each positive integer $k,$ we constructed bolts $l_{k}$ and bolt
functionals $r_{l_{k}}$, for which the inequality (2.32) holds. Note that $%
\left\Vert r_{l_{k}}\right\Vert \leq 1,$ for all $k.$ By the well-known
result of functional analysis (any bounded set in $E^{\ast }$, dual for a
separable Banach space $E$, is precompact in the weak$^{\text{*}}$
topology), the sequence $\left\{ r_{l_{k}}\right\} _{k=1}^{\infty }$ has a
weak$^{\text{*}}$ cluster point. Denote this point by $r_{0}.$ Then from
(2.32) it follows that%
\begin{equation*}
\left\vert r_{0}(h)\right\vert \geq M-\varepsilon .
\end{equation*}

The above inequality completes the proof. \end{proof}

\smallskip

\textbf{Remark 1.} In Theorem 2.3, we use the closedness assumption. For a
given compact space $X,$ the closedness of $A_{1}+A_{2}$ in $C(X)$ strictly
depends on the internal structure of $X.$ There are several results on
closedness of a sum of two algebras (see, e.g., \cite{12,Mar2,Med}). The
most explicit and practical result is due to Medvedev (see \cite{Med}). He
showed that the sum $A_{1}+A_{2}$ is closed in $C(X)$ if and only if the
lengths of all irreducible bolts of $X$ are uniformly bounded. A bolt $%
\{x_{1},...,x_{m}\}$ is called irreducible if there does not exist another
bolt $\{y_{1},...,y_{l}\}$ with $y_{1}=x_{1},$ $y_{l}=x_{m}$ and $l<m$. For
example, the set of functions $\varphi (x)+\psi (y)$ defined on a compact
set $Q\subset \mathbb{R}^{2}$ is closed in the space $C(Q)$ if $Q$ has a
vertical or horizontal bar (a bar is a closed segment in $Q,$ projection of
which into the $x$ or $y$ axis coincides with the projection of $Q$ into the
same axis). This is because any two different points of $Q$ can be connected
by means of a bolt (with respect to the coordinate projections), using at
most $2$ points of a bar. That is, in this case, the lengths of irreducible
bolts are not more than $4$. The similar argument can be applied to a
compact set $Q,$ which contains a broken line $L,$ with segments parallel to
the coordinate axes, such that the projections of $Q$ and $L$ into one of
the coordinate axis coincide. At the same time, there are compact sets of
simple structure, for which the closedness result does not hold. Take, for
example, the triangle $ABC$ with $A=(0,0)$, $B=(2,2)$ and $C=(1,0).$ In this
case, there is no number bounding the lengths of irreducible bolts $%
(A,...,X_{n})\subset ABC,$ $n=1,2,...$, provided that $X_{n}$ tends to $B$,
as $n\rightarrow \infty $.

\bigskip

\textbf{Remark 2.} There is a difference between the convergence $\left\Vert
h_{n}\right\Vert \rightarrow E(h),$ as $n\rightarrow \infty ,$ and the
convergence of the sequence $\{h-h_{n}\}_{n=1}^{\infty }$ to a best
approximation to $h$ from $A_{1}+A_{2}.$ The latter is much stronger.
Theorem 2.3 does not guarantee that the sequence $\{h-h_{n}\}$ converges in
the latter strong sense. We do not know if it converges to a best
approximation to $h$ from $A_{1}+A_{2}.$ The question, for which compact
spaces $X$ and subalgebras $A_{1}$ and $A_{2}$, the sequence of functions
produced by the Diliberto-Straus algorithm converges to a best approximation
is fair, but too difficult to solve in this generality. It should be noted
that Aumann \cite{1} proved that for a rectangle $S=[a,b]\times \lbrack c,d]$
in $\mathbb{R}^{2}$ and for any function $h(x_{1},x_{2})\in C(S),$ the
sequence $\{h-h_{n}\}$ converges uniformly to a best approximation to $h$
from the subspace

\begin{equation*}
D=\left\{ \varphi (x)+\psi (y):\text{ }\varphi \in C[a,b],\text{ }\psi \in
C[c,d]\right\} .
\end{equation*}%
Aumann's proof is based on the equicontinuity of some families of functions,
namely the families $\{u_{n}\}$ and $\{v_{n}\}$ (see the proof of Theorem
2.3). To show the equicontinuity, Aumann substantially uses the
non-expansivness of the averaging operators

\begin{equation*}
(M_{x}f)(y)=\frac{\max_{x}f+\min_{x}f}{2}\text{ and (}M_{y}f)(x)=\frac{%
\max_{y}f+\min_{y}f}{2}.
\end{equation*}

In case of the rectangle $S,$ for any function $f\in C(S)$ and two fixed
points $x_{1},$ $x_{2}\in \lbrack a,b]$ we can always write

\begin{equation*}
\left\vert M_{y}f_{1}-M_{y}f_{2}\right\vert \leq \left\Vert
f_{1}-f_{2}\right\Vert ,\eqno(2.33)
\end{equation*}%
where $f_{1}=f(x_{1},y)$, $f_{2}=f(x_{2},y)$. But inequality (2.33) is not
generally valid for other compact sets $Q\subset \mathbb{R}^{2}$. This is
because for some fixed points $x_{1}$ and $x_{2}$, the projections of the
sets $\{(x_{1},y)\in Q\}$ and $\{(x_{2},y)\in Q\}$ into the $y$ axis do not
coincide; hence $M_{y}f_{1}$ and $M_{y}f_{2}$ are the averaging operators
over two different sets.

\smallskip

\textbf{Remark 3.} Let us explain why we consider the sum of only two
algebras. The problem is that for a sum of more than two algebras, the
appropriate generalization of Diliberto-Straus algorithm does not hold. In
fact, it does not hold even in the simplest case of approximation by sums of
univariate functions. To be more precise, let $E$ be the unit cube in $%
\mathbb{R}^{n}$ and $A_{i}$ be a best approximation operator from the space
of continuous functions $C(E)$ to the subspace of univariate functions $%
G_{i}=\{g_{i}\left( x_{i}\right) :~g_{i}\in C[0,1]\}$, $i=1,...,n.$

That is, for each function $f$ $\in C(X)$, the function $A_{i}f$ is a best
approximation to $f$ from $G_{i}.$ Set
\begin{equation*}
Tf=(I-A_{n})(I-A_{n-1})\cdot \cdot \cdot (I-A_{1})f,
\end{equation*}%
where $I$ is the identity operator. It is clear that
\begin{equation*}
Tf=f-g_{1}-g_{2}-\cdot \cdot \cdot -g_{n},
\end{equation*}%
where $g_{k}$ is a best approximation from $G_{k}$ to the function $%
f-g_{1}-g_{2}-\cdot \cdot \cdot -g_{k-1}$, $k=1,...,r.$ Consider powers of
the operator $T$: $T^{2},T^{3}$ and so on. What is the limit of $\left\Vert
T^{n}f\right\Vert $as $n\rightarrow \infty $ ? One may expect that the
sequence $\{\left\Vert T^{n}f\right\Vert \}_{n=1}^{\infty }$ converges to $%
E(f)$ (the error of approximation from $G_{1}+\cdot \cdot \cdot +G_{n})$, as
in the case $n=2$. This conjecture was first proposed in 1951 by Diliberto
and Straus \cite{4}. But later it was shown by Aumann \cite{1}, and
independently by Medvedev \cite{16}, that the algorithm may not converge to $%
E(f)$ for the case $n>2$.

It should be noted that in Halperin \cite{H} it was proved that the
algorithm converges as desired for any finite number of closed subspaces in
the Hilbert space setting (without demanding closure of the sum of the
subspaces). In fact, in a uniformly convex and uniformly smooth Banach space
it was proved in Pinkus \cite{P} that if the sum is closed then the
algorithm converges (in the strong sense of the previous remark) for any
finite number of closed subspaces.

\bigskip

From Theorem 2.3 and Proposition 2.1 (see also Example 2) one can obtain the
following corollary which is a generalization of the classical Diliberto and
Straus theorem from rectangular sets to special compact convex sets of $%
\mathbb{R}^{2}$.

\bigskip

\textbf{Corollary 2.4.} \textit{Let $Q\subset \mathbb{R}^{2}$ be a convex
compact set with the property that any bolt (with respect to the coordinate
projections) in $Q$ can be made closed by adding only a fixed number of
points of $Q$. Let, in addition, $A_{1}=\{f(x)\}$ and $A_{2}=\{g(y)\}$ be
the algebras of univariate functions, which are continuous on the
projections of $Q$ into the coordinate axes $x$ and $y$, respectively. Then
for a given function $h\in C(Q)$, the sequence $%
h_{1}=h,h_{2n}=h_{2n-1}-Fh_{2n-1},h_{2n+1}=h_{2n}-Gh_{2n},n=1,2,...$,
converges in norm to the error of approximation from $A_{1}+A_{2}$.}

\bigskip

\textbf{Acknowledgment.} The authors are grateful to the referee for
numerous comments and suggestions that improved the original manuscript.

\end{document}